\documentclass{svjour3}

\usepackage{amscd,amsfonts,amssymb,amsmath,latexsym,array,hhline,xcolor,graphicx}

\usepackage{latexsym}
\usepackage{times}

\newcommand\F{\mbox{I\kern-2pt F}}

\newcommand\cE{{\cal E}}

\newcommand\cF{{\cal F}}

\newcommand\cL{{\cal L}}

\newcommand\cO{{\cal O}}

\newcommand\cP{{\cal P}}

\newcommand\e{{\varepsilon}}

\newcommand{\RNumb}[1]{\uppercase\expandafter{\romannumeral #1\relax}}

\def\E{{\bf E}}
\def\P{{\bf P}}

\def\R{{\mathbb R}}
\def\bbr{{\mathbb R}}

\newcommand\fdem{$\Box$}
\newcommand\beq{\begin{equation}}
\newcommand\eeq{\end{equation}}
\newcommand\bea{\begin{eqnarray}}
\newcommand\eea{\end{eqnarray}}
\newcommand\bean{\begin{eqnarray*}}
\newcommand\eean{\end{eqnarray*}}

\begin{document}

\title{Ruin problems with investments on a finite interval:  PIDEs and their viscosity solutions}

\author{ Viktor Antipov  
\and Yuri Kabanov} 
\date{}
\institute{  \at Lomonosov Moscow State University and ``Vega" Institute, Moscow, Russia \\
\email{stayaptichek@gmail.com}.
\and    Lomonosov Moscow State University and   Laboratoire de Math\'ematiques, Universit\'e de
Franche-Comt\'e, Besan\c{c}on, 
France \\
  \email{ykabanov@unif-fcomte.fr}. 
  }

\date{Received: date / Accepted: date}

\titlerunning{Ruin Problems with Investments}

\maketitle

\begin{abstract}
The study deals with the ruin problem when an insurance company 
  invests  its reserve in a risky asset whose the price dynamics is given by a geometric 
  L\'evy process.  Considering the ruin probability as a  of the capital reserve 
   we obtain for it a partial integro-differential equation understood in a viscosity sense and prove a result on the uniqueness  of the viscosity solution for a corresponding boundary 
  value problem.   
\end{abstract}


 
 \section { Introduction} 

In the classical models of collective risk theory it was assumed that   insurance companies keep their reserves in cash or put  it into a bank account, typically, with zero interest rate.  The dynamic of the reserve was described, in the majority of cases by a compound Poisson process or, in more advanced  
frameworks,  by a L\'evy process. The use of the latter is justified by the fact that the majority of claims are small with respect to the total volume, the interarrival times between small claims are small with respect to the planning horizon and the limit theorems can be applied.  
    
In more realistic recent models it is assumed that insurance companies operate in a financial environment and  their  reserves, totally, or partially, are invested into  risky assets having various types of price dynamic used in quantitative finance.   In this paper we consider a rather general model where the price dynamic  and the insurance business activity  are described by two independent  L\'evy processes. In this case the capital reserve evolution solves  a linear 
stochastic equation and  referred to as the generalized Ornstein--Uhlenbeck process driven by the mentioned L\'evy processes. To the moment there is quite a lot of works  devoted  to  ruin problems with investments, mainly on asymptotic behavior  of the exit probability of the capital reserve process  from the positive half-axis.  Usually papers deal with the asymptotic behavior of the ultimate ruin probability when the initial capital tends to infinity, see, e.g., \cite{KP2020} and references therein.   "Ultimate" means: on an infinite time interval.  For very specific cases it is proven that the ultimate ruin probability is smooth and satisfies, in the classical sense, an integro-differential equation, see \cite{KP}, \cite{KPukh}.   
 Though  from the practical point of view the ruin probability with investments  on a finite interval  is more important, we do not know any works on it.     

It the present paper we obtain  a partial  integro-differential equation (PIDE) for more general functionals
of the capital reserve. This equation is understood in the viscosity sense developing further the ideas from the paper \cite{Belk}.  A uniqueness theorem for a boundary value problem is proven.   Note that the viscosity approach allows us to avoid quite a delicate question of smoothness of the studied functions   with respect to the initial capitaL 
  
\section{The model}
In the description of the model we follow \cite{KP2020} and use the notations of stochastic calculus.  

We are given  a filteblue probability space $(\Omega,\cF,{\bf F}=(\cF_r)_{r\ge 0},\P)$ with 
two independent L\'evy processes  $R$ and $P$, having the 
 L\'evy triplets  $(a,\sigma^2,\Pi)$ and $(a_P,\sigma_P^2,\Pi_P)$, and admitting, respectively, the canonical decompositions   
 \beq
\label{Rt}
R_t=at+\sigma W_t+xI_{\{|x|\le 1\}}*(\mu -\nu)_t+xI_{\{|x|>1\}}*\mu_t
\eeq
and 
\beq
\label{Pt}
P_t=a_Pt+\sigma_P W^P_t+xI_{\{|x|\le 1\}}*(\mu^P -\nu^P)_t+xI_{\{|x|> 1\}}*\mu^P_t,  
\eeq
where $\mu$ and $\mu^P$ are jump measures, $\nu(dt,dx)=\Pi(dx)dt$ and $\nu^P(dt,dx)=\Pi_P(dx)dt$, $W$ and $W^P$ are Wiener processes. Let denote by $R^d$ and $P^d$ purely discontinuous martingale components of $R$ and $P$, i.e. the processes 
$$
R^d:=xI_{\{|x|\le 1\}}*(\mu -\nu),  \qquad P^d:=xI_{\{|x|\le 1\}}*(\mu^P -\nu^P).
$$

We consider the process $X=X^{t,u}$, defined on the interval  $[t,\infty[$ by the formula 
\beq
\label{maintint}
X^{t,u}_r=u+\int_{]t,r]}X^{t,u}_{s-}dR_s+P_r-P_t,   
\eeq 
or, in the traditional differential form,  
\beq
\label{maindiff}
dX^{t,u}_r=X^{t,u}_{r-}dR_r+dP_r, \quad X^{t,u}_t=u. 
\eeq
  
We assume the condition  $\Pi(]\infty,-1])=0$ implying that the stochastic (Dol\'eans) exponential $S_t:=\cE (R)$ is strictly positive (a.s.). 

In the actuarial context $X^{t,u}$ is the process describing the evolution of the reserve of an insurance company starting its activity at time $t>0$ with the initial capital $u>0$, having the time horizon $T>t$, and investing the total of its reserve 
in a risky asset with the price process $S$ solving the stochastic equation 
$dS=S_-dR$, $S_0=1$. 
  That is, $S$ is a geometric L\'evy process; 
in the particular case, where $R$ is a Wiener process with drift, $S$ is a geometric Brownian motion, the most popular model for the stock price. 

Recall that in the classical models 
the process $P$, describing the business activity of the insurance company is, usually, a compound Poisson process with drift. The latter can be represented in the form 
$$
P_r=cr+\sum_{i=1}^{N_r}\xi_i, 
$$
where $N$ is a Poisson process independent on the i.i.d. sequence $\xi_i$.   
The case where $c> 0$ and $\xi_i<0$ corresponds to the non-life insurance. The model  $c< 0$ and  $\xi_i>0$  serves to describe the  life insurance business where the company pays the annuity with the intensity $|c|$ and gets benefits  $\xi_i$ or the capital of a venture capital paying salary and selling innovations. 
 Note also that in the modern literature the generalized Ornstein--Uhlenbeck model  suggested by Paulsen in  \cite{Paul-93}  attracted a lot of attention since  a general L\'evy 
process $P$  is viewed as a reasonable approximation  of the claim flow of a large company with many small customers.  

\smallskip 

 Let  $\tau^{t,u}:=\inf \{s>t\colon X^{t,u}_s\le 0\}$. For a measurable function 
 $V:\bbr \to [0,1]$ vanishing on $]-\infty,0]$ we can define the function 
 \beq
 \label{Psi}
 \Psi(t,u)=\Psi(t,u;T):=\E V(X_{\tau^{t,u}}^{t,u})I_{\{\tau^{t,u}<T\}}.
\eeq
   If 
 $V=I_{]0,\infty[}$, the value  $\Psi(t,u)=\P (\tau^{t,u}<T)$ for $u>0$ is the probability of the exit of $X^{t,u}$ from $]0,\infty[$ before $T$. In particular, in  the actuarial terminology, $\Psi(0,u)$ is 
 the ruin probability on the time interval $[0,T]$ for the company starting its activity at time  zero with the initial capital $u$. In the case, where $T=\infty$, and  the $V=I_{]0,\infty[}$
 function $\Psi(0,u;\infty)$, $u>0$, is called {\it ultimate ruin probability}. It is the latter which is studied in the  absolute majority of works on ruins with investments. 
 
\section{Partial integro-differential equation}
Let $\cL:=\cL^d_R+\cL^i_P+\cL^d_R+\cL^i_P$, where
 \bea
\label{LRd} 
 \cL^d_R\Psi(t,u)&:=&\frac 12 \sigma^2u^2 \Psi_{uu}(t,u)+ au\Psi_u(t,u),\\
 \label{LRi}  
  \cL^i_R\Psi(t,u)&:=&\int (\Psi(t,u(1+z))-\Psi(t,u)-\Psi_u(t,u)uzI_{\{|z|\le 1\}})\Pi(dz),\\
 \label{LPd} 
 \cL^d_P\Psi(t,u)&:=&\frac 12 \sigma_P^2 \Psi_{uu}(t,u)+a_P\Psi_u(t,u),\\
 \label{LPi} 
 \cL^i_P\Psi(t,u)&:=&\int (\Psi(t,u+z)-\Psi(t,u)- \Psi_u(t,u)zI_{\{|z|\le 1\}})\Pi_P(dz). 
 \eea
 
 \begin{theorem}
The function $\Psi(t,u)$  on $]0,T[\times ]0,\infty[$ is a viscosity solution of the 
partial  integro-diffe\-rential  equation 
 \beq
 \label{basiceq}
 \Psi_t(t,u)+\cL\Psi(t,u)=0.
\eeq
 \end{theorem}
 
 The above equation  is a symbolical one since we do not know whether the given function $\Psi$ (defined in (\ref{Psi})) has derivatives.  E.g., the ruin probability  it  is not differentiable even in the case of the Cram\'er model ($R=0$, $P$ is a compound Poisson process with drift, the equation (\ref{basiceq}) is of the first order) if the distributions of claims is discrete.  
 The idea of viscosity solutions is to substitute in the equation the function $\Psi$  by a suitably chosen  smooth test functions.  More precisely, one needs to verify at each point of the domain  each element of a suitably defined set of test functions  is a supersolution of (\ref{basiceq}) and for another set it is a subsolution.

Let denote by $C_b^{1,2}(t,u)$ the set of bounded continuous  functions $f:[0,T[ \times \R\to \R$  continuously differentiable in the time variable and two times  continuously differentiable in the phase variable   
  in some neighborhood of the point
$(t,u)\in [0,T[ \times ]0,\infty[ $ and equal to zero for $u\in ]-\infty,0]$.  For $f\in C_b^{1,2}(t,u)$ the values $\cL_Rf(t,u)$ and $\cL_Pf(t,u)$ are well-defined.

For the operator $\cL$ we introduce the following definitions. 

A  function $\Psi:]0,T[ \times \R\to [0,1]$   is called  {\it viscosity supersolution}
 if  for every  point $(t,u)\in ]0,T[ \times]0,\infty[$ and  every   function $f\in C_b^{1,2}(t,u)$ such that
$f(t,u)= \Psi(t,u)$ and $f \le \Psi$
 the inequality $f_t(t,u)+\cL f(t,u)\le 0$ holds.

A  function $\Psi:[0,T[ \times \R\to  [0,1]$  is called {\it viscosity subsolution} 
 if  for every  point \ $(t,u)\in [0,T[ \times]0,\infty[$ and  every   function $f\in C_b^{1,2}(t,u)$ 
 such that
$f(t,u)=\Psi(t,u)$ and $f\ge \Psi$
 the inequality $f_t(t,u)+\cL f(t,u)\ge 0$ holds.

A  function $\Psi:[0,T[ \times \R\to  [0,1]$  is  {\it viscosity solution}  if $\Psi$ is simultaneously a viscosity super- and
subsolution.

At last, a  function $\Psi\in C_b^{1,2}([0,T]\times [0,1])$  is called  {\it classical
supersolution}  if  the inequality  $(\partial/\partial t+\cL)
\Psi \le 0$ holds on on $[0,T[\times ]0,\infty[$. We add the adjective {\it strict}
when $(\partial/\partial t+\cL)
\Psi<0$ on  $[0,T]\times ]0,\infty[$.

\begin{lemma} 
The ruin probability  $\Psi$ as a function of $(t,u)$ is a viscosity supersolution.
\end{lemma} 
{\sl Proof.} Let $f\in C_b^{1,2}(t,u)$ be such that
$f(t,u)=\Psi(t,u)$ and $f\le \Psi $ on $[0,T]\times \R$. 
Take  $\epsilon\in ]0,u/2[$ and $h\in ]0,T-t[$ such that that the rectangular 
$[t,t+h]\times [u-\epsilon,u+\epsilon]$ is inside a ball  in which the function $f$ is in  $C^{1,2}$. 
  
Let $M=(M_s)_{s\ge t}$ be the local martingale with 
$$
dM_s=dX^c_s+X_{s-}dR^d_s+dP^d_s, \qquad M_t=0,
$$  
where $X^c$ is the continuous martingale component of $X$,  
$$
dX^c_s:=\sigma X_{s-}dW_s+\sigma_PdW^P_s,  \qquad d\langle X^c\rangle_s= \sigma^2 X^2_{s}ds+\sigma_P^2ds. 
$$

Define the stopping time $\tau_h=\tau^{\epsilon}_h\le t+h$ by putting 
$$
\tau_h:=\inf\big\{s> t\colon  X^{t,u}_{s}
\notin[u-\epsilon\,,u+\epsilon]\big\}\wedge \inf\{s> t\colon  |M_{s}|\ge 1\}\wedge (t+ h). 
$$
By  the Ito formula in its standard form applied to the semimartingale $X=X^{t,u}$
\bean
 f(\tau_h,X_{\tau_h})&=&f(t,u)+\int_{]t,\tau_h]}f_u(s,X_{s-})dX_s+\int_{]t,\tau_h]}f_t(s,X_{s-})ds  
 \\
 &&
 +\frac 12 \int_{]t,\tau_h]}f_{uu}(s,X_{s-})d\langle X^c_s\rangle\\
 &&+
 \sum_{s\in ]t, \tau_h]}\big(f(s,X_{s-}+\Delta X_{s})-f(s,X_{s-})-f_u(s,X_{s-})\Delta X_s\big).
\eean

Representing  the
integral with respect to $dX$   as  
$$
\int_{]t,\tau_h]}f_u(s,X_{s-})dM_s+
\int_{]t,\tau_h]}f_u(s,X_{s-})(aX_{s-}+a_P)ds+ J^R+J^P,
$$
where
\bean
J^R&:=&\int_{]t,\tau_h]}\int f_u(s,X_{s-})X_{s-}zI_{\{|z|>1\}}\mu(ds,dz), \\
J^P&:=&\int_{]t,\tau_h]}\int f_u(s,X_{s-})zI_{\{|z|>1\}}\mu^P(ds,dz). 
\eean
and  the sum of jumps as
\bean
&&\int_{]t,\tau_h]}\int \big(f(s,X_{s-}(1+z))-f(s,X_{s-})-f_u(s,X_{s-})X_{s-} zI_{\{|z|\le 1\}}\big)\mu(ds,dz)\\
&&
+\int_{]t,\tau_h]}\int \big(f(s,X_{s-}+z)-f(s,X_{s-})-f_u(s,X_{s-})zI_{\{|z|\le 1\}}\big)\mu^P(ds,dx)-J^R-J^P, 
\eean 
we rewrite  the right-hand side of the Ito formula  in the form 
\bean
 &&f(t,u)+\int_{]t,\tau_h]}f_u(s,X_{s-})(dX^c_s+X_{s-}dR^d+dP^d)\\
 &&+\int_{]t,\tau_h]}\Big(f_t(s,X_{s-})+f_u(s,X_{s-})(aX_{s-}+a_P)
+\frac 12f_{uu}(s,X_{s-})(\sigma^2 X^2_{s-}+\sigma_P^2)\Big) ds\\
 && + \int_{]t,\tau_h]}\int \big(f(s,X_{s-}(1+z))-f(s,X_{s-})-f_u(s,X_{s-})X_{s-} zI_{\{|z|\le 1\}}\big)\mu(ds,dz)\\
 &&+\int_{]t,\tau_h]}\int \big(f(s,X_{s-}+z)-f(s,X_{s-})-f_u(s,X_{s-})zI_{\{|z|\le 1\}}\big)\mu^P(ds,dz)
\eean
Taking the expectation we get that 
$$
\E f(\tau_h,X_{\tau_h})=f(t,u)+\E \int_t^{\tau_h}\big (f_t(s,X_s)+ \cL(s,X_s))ds. 
$$ 
Recall that $f(t,u)=\Psi(t,u)$ and $f\le \Psi$. The crucial observation is that 
$$
\Psi(t,u)=\E\,\Psi(\tau_h,X^{t,u}_{\tau_h})
$$ 
due to the strong Markov property of $X$. Combining these facts, we get that 
$$
\E \int_{t}^{\tau_h}\big (f_t(s,X_{s})+\cL
f(s,X_{s})\big)ds\le 0.
$$
Almost surely  $\tau_h(\omega)=t+h$ when $h$ is sufficiently small  
(the threshold  for which this equality holds depends on $\omega$). 
Using the Lebesgue theorem on dominated convergence we get that 
$$
\cL f(t,u)=\lim_{h\downarrow 0}\frac 1h \E \int_{t}^{\tau_h}\big (f_t(s,X_{s})+\cL
f(s,X_{s})\big)ds\le 0
$$
and the lemma is proven. \fdem

\begin{lemma} 
The function $\Psi$ is a viscosity subsolution.
\end{lemma} 
{\sl Proof.}  One argue as above but with the inequality $f\ge \Psi$. \fdem

\begin{remark} If the function $\Psi\in C^{1,2}$, then the equation $\cL\Psi(t,u)=0$
holds in the classical sense. 
\end{remark}

\section{Parabolic jets}
Let $f$ and $g$ be functions defined in a neighborhood of the origin of $\R^2$. Denoting $(h_1,h_2)$ a generic point of $\R^2$, we
shall write $f(.)\lessapprox g(.)$ if $f(h_1,h_2)\le g(h_1,h_2)+o(|h_1|+|h_2|^2)$ as
$|h|\to 0$. 
The notations $f(.)\gtrapprox g(.)$ and $f(.)\approx
g(.)$ have the obvious meaning. 

For $(b,p,A)\in \R^2\times \R_+$  we consider the 
function
$$
Q_{b,p,A}(t,x):= bt + px+(1/2)Ax^2,\quad (t,x)\in {\R}^2,
$$
linear in $t$ and quadratic in $x$.  Define  the parabolic {\it super}- and {\it subjets} of a function $v$ 
at the point $(t,x)\in K:= ]0,T[\times ]0,\infty[$: 
\bean
\cP^+v(t,x)&:=&\{(b,p,A)\in \R^3:\ v((t,x)+.)\lessapprox v(t,x)+Q_{b,p,A}(.)\},\\
\cP^-v(t,x)&:=&\{(b,p,A))\in \R^3:\ v((t,x)+.)\gtrapprox v(t,x)+Q_{b,p,A}(.) \}. 
\eean

The following lemma asserts that for $v\in C(K)$  and any element  $(b,p,A)\in \cP^+v(t,x)$,  
$(t,x)\in K$, 
one can associate a test function $f$ dominating $v$ on $K$ and uniformly close to $v$,  having at the point $(t,x)$ the same value as $v$, belonging to  $C^{1,2}$ in a neighborhood of  $(t,x)$ and having at this point the first derivative in $t$ equal to $b$ and the first and second derivatives in $x$, equal, respectively, to $p$ and $A$. Moreover, $f=v$ outside any arbitrary small neighborhood of $(t,x)$.

\begin{lemma}
 \label{jet1}
 Let $v\in C(K)$, $(t,x) \in  K$, and   $(b,p,A)\in \cP^+v(t,x)$.  Then for  $\e,\delta\in ]0,1]$  there exist a real  $r>0$ depending on $\e$ and  a function 
 $f\in C^{1,2}(t,x)$ depending on $\e$ and $\delta$  with the following properties: 
 
 \smallskip
 $(i)$ $f(t,x)=v(t,x)$, $f_t(t,x)=b$, $f_x(t,x)=p$,  $f_{xx}(t,x)=A$; 
 
  \smallskip
 $(ii)$ $v\le f\le v+\e$ on $K$;
 
  \smallskip
 $(iii)$ $f((t,x)+h)=v(t,x)+Q_{b,p,A}(h)+\rho(h)$ for $h\in \cO_{r}(0,0)$, where $C^{1,2}$-function $\rho: {\R}^2\to {\R}$ with compact support is such that 
\beq
\label{jet}
 \lim_{|h|\to 0}\frac{\rho(h)} {|h_1|+|h_2|^2}=0,    
\eeq
 \smallskip
 $(iv)$ $f=v$ on the set $K\setminus \cO_{r'}(t,x)$ where $r':=r+\delta$. 
\end{lemma}

\noindent
{\sl Proof.} Take $r\in ]0,1[$ such that the ball $\cO_{2r}(t,x)\subset K$.
By definition of  $\cP^+v(t,x)$, 
$$
\varrho (h):=v((t,x)+h)-v(t,x)-Q_{b,p,A}(h)\le (|h_1|+|h_2|^2) \varphi (h),
$$
where $\varphi(h)\to 0$ as $h\to  0$. We consider on $\cO_{2r}(0,0)\cap \R_+^2$
the  function
$$
\widehat\varphi(u):=\sup_{\{h:\ |h_1|\le u_1,|h_2|\le u_2\}}\frac {\varrho^+(h)}{|h_1|+|h_2|^2}
\le \sup_{\{h:\ |h_1|\le u_1,|h_2|\le u_2\}}\varphi^+(u). 
$$
Obviously, $\widehat\varphi $ is continuous,  increasing in $u_1$, $u_2$, and $\widehat\varphi(u)\to
0$ as $u\to  0$. We extend $\widehat\varphi $ to a positive continuous function on ${\R}^2_+$
with $\widehat\varphi(u)=0$ for $|u|\ge 1$ and associate with it the function $\Upsilon:\R^2_+\to \R_+ $ by putting  
$$
\Upsilon (u)=\Upsilon(u_1,u_2):=\int_{u_1}^{2u_1}\widehat\varphi(\xi,u_2)d\xi  +   \frac 2{3}\int _{u_2}^{2u_2}\int _\eta^{2\eta}\widehat\varphi (u_1,\xi)d\xi
d\eta. 
$$
Then 
$$
(u_1+u_2^2) \widehat\varphi (u)\le \Upsilon(u)\le (u_1+u_2^2) \widehat\varphi (4 u). 
$$
It follows that the function 
$\rho: h\mapsto \Upsilon(|h_1|,|h_2|)$ has 
a compact support, belongs to $C^{1,2}(\cO_{r}(0,0))$, it vanishes at the origin together with the partial derivatives of the first order in $h_1$ and $h_2$  and the second partial derivative in $h_2$. For $h\in \cO_{r}(0,0)$ we have that 
$$
v((t,x)+h)-v(t,x)-Q_{b,p,A}(h)\le (|h_1|+|h_2|^2)\widehat \varphi (h)\le \Upsilon (|h_1|,|h_2|)=\rho(h).
$$
Thus, the function 
$$
f^0: (s,y)\mapsto v(t,x)+Q_{b,p,A}(s-t,y-x)+\rho (s-t,y-x)
$$ 
dominates 
$v$ on the ball $\cO_{r}(t,x)$. Without loss of generality, diminishing $r$ if necessary, we may assume that $f^0\le v+\e$ on the ball  $\cO_{r'}(t,x)$ for some $r'>r$.

Let us consider a continuous function $\xi_{r,r'}: {\R}^2\to [0,1]$ such that  $\xi_{r,r'}=1$
on $\cO_r(0,0)$ and $\xi_{r,r'}=0$ outside  $\cO_{r'}(0,0)$. 
The function  
$$
f: (s,y)\mapsto f^0(s,y)\xi_{r,r'}(s-t,y-x)+v(s,y)\big(1-\xi_{r,r'}(s-t,y-x)\big)
$$
has all the needed properties. \fdem

\smallskip
The corresponding assertion for $\cP^-v(x)$  holds  with obvious changes in the formulation.  

The following result is a scalar version of the Ishii lemma. 
\begin{lemma}
\label{Ishii} Let $\cO$ be  an open subset of $\R_+$ and 
$v$ and $\tilde v$ be two continuous functions
on $ ]0,T[\times \cO$. Let  $\Delta_n (t,x,y):= v(t,x)-\tilde
v(t,y)-\frac 12 n|x-y|^2$, where the parameter $n>0$. Suppose that the function $\Delta_n$ attains
a  
maximum at $(t_n, x_n, y_n)\in ]0,T[\times \cO\times \cO$. Then there are reals $b_n$, $A_n$, and $\tilde A_n$ 
such that
$$
(b_n,n(x_n-y_n),A_n)\in   \bar \cP^+v( t_n, x_n), \qquad (b_n,n( x_n- y_n),\tilde A_n)\in   \bar \cP^-\tilde v(t_n, y_n),
$$
and, for all $x,y\in \cO$ 
\beq
\label{doublematrix}
A_n x^2- \tilde A_n y^2\le 3n(x-y)^2. 
\eeq
In particular, $A_n\le \tilde A_n$. 
\end{lemma}

\begin{remark} Recall that the symbol 
 $(b_n,n(x_n-y_n),A_n)\in   \bar \cP^+v( t_n, x_n)$ means that there are sequences 
$(t_{nk},x_{kn})$ and $(b_{nk},p_{nk},A_{nk})\in  \cP^+v( t_{n,k}, x_{nk})$, $k\ge 1$, such that 
$$
(t_{nk},x_{kn})\to (t_{n},x_{n}), \quad (b_{nk},p_{nk},A_{nk})\to (b_n,n(x_n-y_n),A_n), \ \ \ k\to \infty.
$$
Similar meaning has the symbol $(b_n,n( x_n- y_n),\tilde A_n)\in   \bar \cP^-\tilde v(t_n, y_n)$. 
\end{remark}

\section{The maximum principle and uniqueness theorem } 
\label{Maximum} 
\begin{theorem}
\label{maximum}
Let   $\psi ,\tilde \psi$ be continuous functions on $[0,T]\times \bar R_+$ with values in the interval $[0,1]$ be, respectively,  viscosity subsolution and supersolution of   (\ref{basiceq}) such that  
\beq
\label{BC}
\psi(t,0) = \tilde\psi(t,0), \quad \psi(t,\infty)=\tilde\psi(t,\infty), \quad  \psi(T,u) = \tilde\psi(T,u).
\eeq
Then $\tilde \psi\ge  \psi$. 
\end{theorem}
\noindent {\sl Proof.} 
For any $\delta>0$ the function $\tilde \psi^\delta(t,u):= \tilde \psi(t,u)+\delta/t\ge \delta/t$ is  a strict supersolution of (\ref{basiceq}). It converges uniformly to infinity as $t\to 0$.   We prove that $\tilde \psi^\delta \ge  \psi$ and get the result by letting $\delta$ tends to zero.    

If the claimed inequality fails, then 
 $\psi(t_0,u_0)=\tilde \psi^\delta(t_0,u_0)+\e_0$ for some $u_0\in ]0,\infty[$, $t_0\in [0,T[$, and $\e_0>0$. 
 
 \smallskip

{\bf 1.} The first step in a way to get a contradiction consists to use the method of doubling variables and  apply the Ishii lemma to each element of the sequence of functions $\Delta _n=\Delta^\delta_n$ attaining their  maximum values.  

Define  a sequence of continuous functions $\Delta_n:  ]0,T]\times \R_+^2 \to \R$, 
$n\ge 0$, by 
$$
\Delta_n(t,u,w):=\psi(t,u)-\tilde \psi^\delta(t,w)-\frac 12 n |u-w|^2. 
$$
Note that $\Delta_n(t,u,u)=\Delta_0(t,u,u)$ for all $(t,u)\in ]0,T]\times\R_+$ and
$\Delta_0(t,0,0)\le 0$. 

Let us consider the non-empty closed sets 
$$
\Gamma_n:=\{(t,u,w)\in ]0,T]\times \R^2_+ \colon \Delta_n(t,u,w)\ge \e_0\}. 
$$ 
Then $\Gamma_{n+1}\subseteq \Gamma_n \subseteq [\delta,T]\times \R^2_+ $.   The set $\Gamma_1$ is bounded. Indeed, if it is not the case, we can find a sequence $(s_m,u_m,w_m)\in \Gamma_1$ such that $s_m\to s$,  $u_m,w_m\to \infty$.
Due to assumed continuity of $\psi$ and 
$\tilde \psi$ we get the inequality
$$
\limsup_m \Delta_1(s_m,u_m,w_m)\le -\delta/s <0,
$$
i.e.  a contradiction.

It follows that 
all functions $\Delta_n$ attain their  maximal values over $[0,T]\times \R^2_+$ on the compact subset $\Gamma_1$. 
That is, there exist
$(t_n,x_n,y_n)\in[\delta,T]\times \R^2_+ $ such that
\bean
\Delta_n(t_n,x_n,y_n)&=&\bar \Delta_n:=  \sup_{(t,x,y)\in [0,T]\times \R^2_+ }\Delta_n(t,x,y)\\
&\ge& \bar \Delta
  :=\sup_{(t,x)\in [0,T]\times \R_+ }\Delta_0(t,x,x)\ge \e_0>0.
\eean
This implies that the sequence $n|x_n-y_n|^2$ is bounded.
We continue to argue (without introducing new
notations) with a subsequence of indices along which $(t_n,x_n,y_n)$ converge to
some limit $(\widehat t,\widehat x,\widehat x)$. Necessarily,
$n|x_n- y_n|^2\to 0$ (otherwise
 $\Delta_0 (\widehat t,\widehat x,\widehat x)>\bar \Delta$).
 Also, it is easily seen that
$\bar \Delta_n \downarrow \Delta_0(\widehat t,\widehat x,\widehat x)=\bar \Delta$.
Thus,
 $(\widehat t,\widehat x)\in  K:=]0,T[\times ]0,\infty[$.  
\smallskip

{\bf 2.} The second step consists in constructing test functions. 

By virtue of the Ishii lemma  applied to the functions $\psi 
$ and $\tilde \psi^\delta $ at the point $(t_n,x_n, y_n)$ there
are  numbers $b_n$, $A_n$ and $\tilde A_n$ satisfying (\ref{doublematrix})  such
that
\beq
\label{natural}
(b_n,n(x_n-y_n),A_n)\in   \bar \cP^+\psi(t_n,x_n), \qquad (b_n,n(x_n-y_n),\tilde A_n)\in   \bar \cP^-\tilde \psi^\delta (t_n, y_n).
\eeq

 According to the remark following  the Ishii lemma, there are  points 
$(t_{nk},x_{nk})\in K$ and $(b_{nk},p_{nk},A_{nk})\in  \cP^+\psi( t_{nk}, x_{nk})$ such that  
$$
(t_{nk},x_{kn})\to (t_{n},x_{n}), \quad (b_{nk},p_{nk},A_{nk})\to (b_n,n(x_n-y_n),A_n), 
$$
and $(t'_{nk},y_{kn})\in K$, $(b'_{nk},p'_{nk},\tilde A'_{nk})\in  \cP^-\tilde \psi^\delta ( t'_{nk}, y_{nk})$  such that  
$$
(t'_{nk},y_{kn})\to (t_{n},y_{n}), \quad (b'_{nk},p'_{nk},\tilde A'_{nk})\to (b_n,n(x_n-y_n),\tilde A_n), \quad k\to \infty.
$$
 
From the elements of  jets $\cP^+\psi(t_{nk},x_{nk})$ and $\cP^-\tilde \psi^\delta(t'_{nk},y_{nk})$
we built functions $f^{nk}$ and $\tilde f^{nk}$ to test, respectively, that  $\psi$ is a viscosity  subsolution of the equation (\ref{basiceq}) at the point $(t_{nk},x_{nk})$ and $\tilde \psi^\delta$ is a viscosity  supersolution   at the point $(t'_{nk},y_{nk})$. 

Take arbitrary positive $\e_{n}, \delta_n \to 0$ as $n \to \infty$ and $\e_{nk} \downarrow \e_{n}, \delta_{nk} \downarrow \delta_n$ as $k \to \infty$ (the particular choice of parameters will be further). By virtue of  Lemma \ref{jet1} 
there are $r_{nk}>0$ and  $f^{nk}\in C_b^{1,2}(t_{nk},x_{nk})$  with the properties  $\psi \le f^{nk}\le \psi + \e_{nk}$,  
\bean
f^{nk}(t_{nk},x_{nk})&=&\psi (t_{nk},x_{nk}), \quad f^n_t(t_{nk},x_{nk})=b_{nk},\quad f^{nk}_x(t_{nk},x_{nk})=p_{nk},\\
&&f^{nk}_{xx}(t_{nk},x_{nk})=A_{nk}, 
\eean
$f^{nk}((t_{nk},x_{nk})+h)=\psi (t_{nk},x_{nk})+Q_{b_{nk},p_{nk},A_{nk}}(h)+\rho_{nk}(h)$ for $h\in \cO_{r^{nk}}(0,0)$, where $C^{1,2}$-function $\rho_{nk}: {\R}^2\to {\R}$ with compact support is such that 
\beq
\label{jeta}
 \lim_{|h|\to 0}\frac{\rho_{nk}(h)} {|h_1|+|h_2|^2}=0,    
\eeq
 \smallskip
 and  $f^{nk}=\psi$ on the set $K\setminus \cO_{r'}(t_{nk},x_{nk})$ where $r'_{nk}:=r_{nk}+\delta_{n}$. 

Similarly, there are $\tilde r_{nk}>0$
and  a function $\tilde f^{nk}\in C_b^{1,2}(t'_{nk}, y_{nk})$  with the properties  $\tilde \psi^\delta \ge \tilde f^{nk}\ge \tilde \psi^\delta -\e_{nk}$,  
\bean
\tilde f^{nk}(t'_{nk}, y_{nk})&=&\tilde \psi^\delta (t'_{nk}, y_{nk}), \quad \tilde f^{nk}_t(t'_{nk},y_{nk})=b'_{nk}, \quad \tilde f^{nk}_x(t'_{nk},y_{nk})=p'_{nk},\\
&& \tilde f^{nk}_{xx}(t'_{nk},y_{nk})=\tilde A'_{nk}, 
\eean
$\tilde f^{nk}((t'_{nk},y_{nk})+h)=\tilde \psi (t'_{nk},y_{nk})+Q_{b'_{nk},p'_{nk},\tilde A'_{nk}}(h)+\tilde \rho_{nk}(h)$ for $h\in \cO_{\tilde r^{nk}}(0,0)$, where $C^{1,2}$-function $\tilde \rho_{nk}: {\R}^2\to {\R}$ with compact support is such that 
\beq
\label{jetb}
 \lim_{|h|\to 0}\frac{\tilde \rho_{nk}(h)} {|h_1|+|h_2|^2}=0,    
\eeq
and $\tilde f^{nk}=\tilde \psi^\delta$ on the set $K\setminus \cO_{\tilde r'_{nk}}(t'_{nk},y_{nk})$ where $\tilde r'_{nk}:=\tilde r_{nk}+\tilde \delta_{n}$.

\smallskip
By definitions of sub- and supersolutions we have that  
\bean
f^{nk}_{t}(t_{nk},x_{nk}) + \cL f^{nk}(t_{nk}, x_{nk})&\ge& 0\\
  \tilde f^{nk}_{t}(t'_{nk},y_{nk}) + \cL \tilde f^{nk}(t'_{nk}, y_{nk})&\le& -\delta/T^2 .
\eean
Hence, 
\beq
\label{bnk}
B_{nk}:= f^{nk}_{t}(t_{nk},x_{nk}) -  \tilde f^{nk}_{t}(t'_{nk},y_{nk}) +   \cL f^{nk}(t_{nk}, x_{nk})-\cL \tilde f^{nk}(t'_{nk}, y_{nk})\ge \delta/T^2 .
\eeq

{\bf 3.} The third step consists in getting a contradiction with our hypothesis. 

After regrouping terms we get that $ B_{nk}= B^{(1)}_{nk}+ B^{(2)}_{nk}$ where 
\bean 
 B^{(1)}_{nk}&:=&\frac 12 \sigma^2 (A_{nk}x_{nk}^2 -\tilde  A'_{nk}y_{nk}^2 )+\frac 12 \sigma^2_P (A_{nk} -\tilde  A'_{nk})
 +a(p_{nk}x_{nk}-p'_{nk}y_{nk})\\
 &&+a_P(p_{nk}-p'_{nk})
+b_{nk}-b'_{nk},\\
B^{(2)}_{nk}
&:=&   \cL^i_Rf^{nk}(t_{nk},x_{nk}) - \cL^i_R\tilde f^{nk}(t'_{nk}, y_n) 
 + \cL^i_Pf^{nk}(t_{nk},x_{nk}) - \cL^i_P\tilde f^{nk}(t'_{nk},y_{nk}) 
\eean
with  $\cL^i_R$ and $\cL^i_P$  defined in (\ref{LRi}) and (\ref{LPi}). 

Passing to the limit as $k\to \infty$ we get that 
$$
\lim_{k} B^{(1)}_{nk}=\frac 12 \sigma^2 (A_nx_n^2-\tilde A_ny_n^2)+\frac 12 \sigma^2_P (A_n-\tilde A_n)+ a n(x_n-y_n)^2. 
$$
By virtue of the Ishii lemma $A_n-\tilde A_n\le 0$,  the first term is dominated by a constant multiplied by $n|x_n-y_n|^2$
and, therefore, converges to zero as $n\to \infty$.  It follows that $\limsup_n\lim_kB^{(1)}_{nk}\le 0$. 

To get a contradiction with (\ref{bnk}) it is sufficient to show that
\bea  
\label{J1} 
&&\limsup_n\limsup_k\big( \cL^i_Rf^{nk}(t_{nk}, x_{nk}) - \cL^i_R\tilde f^{nk}(t'_{nk}, y_{nk}) \big)\leq 0,\\
\label{I1} 
&&\limsup_n\limsup_k\big( \cL^i_Pf^{nk}(t_{nk}, x_{nk}) - \cL^i_P\tilde f^{nk}(t'_{nk}, y_{nk}) \big)\leq 0.
\eea

Let  
\bean 
 H_{nk}(z)&:=&[f^{nk}(t_{nk}, x_{nk}+x_{nk}z)-f^{nk}(t_{nk},x_{nk})-f^{nk}_x(t_{nk},x_{nk})x_{nk}zI_{\{|z|\le 1\}}]\\
&&-[\tilde f^{nk}(t'_{nk},y_{nk}+y_{nk}z)-\tilde f^{nk}(t'_{nk},y_{nk})-\tilde f^{nk}_x(t'_{nk},y_{nk})y_{nk}zI_{\{|z|\le 1\}}], 
\\
G_{nk}(z)&:=&[f^{nk}(t_{nk}, x_{nk}+z)-f^{nk}(t_{nk},x_{nk})-f^{nk}_x(t_{nk},x_{nk})zI_{\{|z|\le 1\}}]\\
&&-[\tilde f^{nk}(t'_{nk},y_{nk}+z)-\tilde f^{nk}(t'_{nk},y_{nk})-\tilde f^{nk}_x(t'_{nk},y_{nk})zI_{\{|z|\le 1\}}].
\eean

In this notation the claimed relations (\ref{J1}) and (\ref{I1}) have the form
\bea
&&\limsup_n\limsup_k\int  H_{nk}(z)\Pi(dz)\le 0,
 \\
&&\limsup_n\limsup_k\int  G_{nk}(z)\Pi_P(dz)\le 0.
\eea 

 In the sequel we  work with $n$ and $k$ large enough  to ensure that $\kappa \le x_{nk}, y_{kn}\le 1/\kappa $ for some (small) constant $\kappa >0$.  
 
On the set  $\{\vert z \vert > 1\}$  the functions  $H_{nk}$  and $G_{nk}$ admit a uniform bounds 
\bean
H_{nk}(z)&\le& \psi(t_{nk}, x_{nk}+x_{nk}z)-\tilde \psi^{\delta}(t'_{nk}, y_{nk}+y_{nk}z)-\psi(t_{nk}, x_{nk})+\tilde \psi^{\delta}(t'_{nk}, y_{nk})+2\e_{nk}, \\ \\
G_{nk}(z)&\le& \psi(t_{nk}, x_{nk}+z)-\tilde \psi^{\delta}(t'_{nk}, y_{nk}+z)-\psi(t_{nk}, x_{nk})+\tilde \psi^{\delta}(t'_{nk}, y_{nk})+2\e_{nk}.
\eean
Taking $\e_{nk}=1/n+1/k$ we get that 
\bean
\limsup_k H_{nk}(z)&\le &\psi(t_{n}, x_{n}+x_{n}z)-\tilde \psi^{\delta}(t_{n}, y_{n}+y_{n}z)-\psi(t_{n}, x_{n})+\tilde \psi^{\delta}(t_{n}, y_{n}) + 1/n. 
\eean

Notice that on the set $\{\vert z \vert > 1\}$ 
\bean
\psi(t_{n}, x_{n}+x_{n}z)-\tilde \psi^{\delta}(t_{n}, y_{n}+y_{n}z) &=& \Delta_n( t_n,x_{n}+x_{n}z,y_{n}+y_{n}z) + (1/2) nz(z + 2) |x_n - y_n|^2 \\
&& + (1/2) n |x_n - y_n|^2 \\
&\le& \Delta_n( t_n,x_{n}+x_{n}z,y_{n}+y_{n}z)  + (3/2) nz^2 |x_n - y_n|^2 \\
&& + (1/2) n |x_n - y_n|^2.
\eean
Hence, on the set $\{\vert z \vert > 1\}$ 
\bean
\limsup_k H_{nk}(z) &\le &\Delta_n( t_n,x_{n}+x_{n}z,y_{n}+y_{n}z)-\Delta_n( t_n,x_{n},y_{n})
+ \frac32 n z^2 |x_n-y_n|^2 + \frac1n. 
\eean
The first term in the inequality above is dominated by the second term and the last two terms converge to zero as $n \to \infty$. 

Similarly, we get that 
\bean
\limsup_k G_{nk}(z)&\le &\psi(t_{n}, x_{n}+z)-\tilde \psi^{\delta}(t_{n}, y_{n}+z)-\psi(t_{n}, x_{n})+\tilde \psi^{\delta}(t_{n}, y_{n})\\
&\le &\Delta_n( t_n,x_{n}+z,y_{n}+z)-\Delta_n( t_n,x_{n},y_{n}) \le 0.
\eean
Using the 	boundedness of $H_{nk}(z) I_{\{\vert z \vert > 1\}}$ and $G_{nk}(z) I_{\{\vert z \vert > 1\}}$ and applying the Fatou lemma we get the bounds
\bea
&&\limsup_n\limsup_k \int H_{nk}(z)I_{\{\vert z \vert > 1\}}\Pi(dz)\le 0,
 \\
&& \limsup_n \limsup_k \int G_{nk}(z)I_{\{\vert z \vert > 1\}}\Pi_P(dz)\le 0.
\eea
 
Take $r_{nk}\in ]0,\kappa [$ small enough to guaranty that for $\theta\in ]0,r_{nk}]$ we have the representations 
\bean
f^{nk}(t_{nk},x_{nk}+\theta)&=&\psi (t_{nk},x_{nk})
+p_{nk}\theta
+ (1/2)A_{nk}\theta^2 +R_{nk}(\theta), \\
\tilde f^{nk}(t'_{nk},y_{nk}+\theta)&=&\tilde \psi^\delta (t'_{nk},y_{nk})
+p'_{nk}\theta
+ (1/2)\tilde A'_{nk}\theta^2 + \tilde R_{nk}(\theta), 
\eean 
where $R_{nk}(\theta)=\rho_{nk}(0,\theta)$,  $\tilde R_{nk}(\theta)=\tilde \rho_{nk}(0,\theta)$.

On the set $U_{nk} := \{z\colon |z|\le \kappa r_{nk}\}\cap \{z\colon |z|\le 1\}$
\bean
H_{nk}(z)&=&
(1/2)(A_{n} x^2_{nk} -\tilde A_{n} y^2_{nk})z^2\\
&& +(1/2)\Big((A_{nk} -A_n)x^2_{nk} - (\tilde A'_{nk} -\tilde A_n)y^2_{nk}\Big)z^2 +R_{nk}(z{x_{nk}})-\tilde R_{nk}(z{y_{nk}}), 
\\
G_{nk}(z)&=&
(1/2)(A_{nk} - \tilde A'_{nk})z^2 + R_{nk}(z)-\tilde R_{nk}(z).
\eean 
and we have the bounds
\bean
H_{nk}(z)&\le &(3/2)n(x_{nk}-y_{nk})^2z^2+ (1/(2\kappa^2))\big(|A_{nk} -A_n|+|\tilde A'_{nk} -\tilde A_n|\big)z^2\\
&&+|R_{nk}(z{x_{nk}})|+|\tilde R_{nk}(z{y_{nk}})|, 
\\
G_{nk}(z)&\le &(1/2) |A_{nk} - \tilde A'_{nk}| z^2 + |R_{nk}(z)|+|\tilde R_{nk}(z)|.
\eean

The ratios $R_{nk}(zx_{nk})/z^2$, $\tilde R_{nk}(zy_{nk})/z^2$, $R_{nk}(z)/z^2$ and $\tilde R_{nk}(z)/z^2$ are bounded and tend to zero as $z\to 0$. 
Thus, for all sufficiently large $n,k$ 
\bean
&& \int_{
U_{nk}} H_{nk}(z)\Pi(dz)\le \frac 1n,
 \\
&& \int_{
U_{nk}} G_{nk}(z)\Pi_{P}(dz)\le \frac 1n. 
\eean

Also we can find $\delta_n > 0$ small enough and $r'_{nk} := \kappa r_{nk} + \delta_n$ (from Lemma \ref{jet1}) such that for sufficiently large $k, n$ on the set $V_{nk} := \{z\colon \kappa r_{nk} \le |z|\le r'_{nk}\} \cap \{z\colon |z|\le 1\}$ we get

\bean	
&& \int_{
 V_{nk}} H_{nk}(z)\Pi(dz)\le \frac 1n,
 \\
&& \int_{V_{nk}} G_{nk}(z)\Pi_{P}(dz)\le \frac 1n.  
\eean

On the set $\{z\colon |z| \ge {r'_{nk}}\}$ 
$$
f^{nk}(t_{nk},x_{nk}+zx_{nk})=\psi(t_{nk},x_{nk}+zx_{nk}),     \quad 
\tilde f^{nk}(t'_{nk},y_{nk}+zy_{nk})=\tilde \psi^\delta(t'_{nk},y_{nk}+zy_{nk}),
$$
$$
f^{nk}(t_{nk},x_{nk}+z)=\psi(t_{nk},x_{nk}+z),     \quad 
\tilde f^{nk}(t'_{nk},y_{nk}+z)=\tilde \psi^\delta(t'_{nk},y_{nk}+z).
$$
implying that  on the set $Y_{nk} := \{z\colon  r'_{nk}\le |z| \le 1\}$

\bean
H_{nk}(z)&=& \psi(t_{nk}, x_{nk}+x_{nk}z)-\tilde \psi^{\delta}(t'_{nk}, y_{nk}+y_{nk}z)-\psi(t_{nk}, x_{nk})+\tilde \psi^{\delta}(t'_{nk}, y_{nk})\\
&&-  (p_{nk}x_{nk}-p'_{nk}y_{nk})z,\\
G_{nk}(z)&=& \psi(t_{nk}, x_{nk}+z)-\tilde \psi^{\delta}(t'_{nk}, y_{nk}+z)-\psi(t_{nk}, x_{nk})+\tilde \psi^{\delta}(t'_{nk}, y_{nk}) \\
&&-  (p_{nk} -p'_{nk})z.
\eean 
The functions $H_{nk}(z)$, $G_{nk}(z)$ are continuous on the compact set $\{z\colon  \delta_{n}\le |z| \le 1\} \supseteq Y_{nk}$, and hence uniformly continuous. Together with the facts that 
$$p_{nk}x_{nk}-p'_{nk}y_{nk}\to n(x_n-y_n)^2, \quad p_{nk} -p'_{nk} \to  0, \quad k\to \infty,$$ we get that for sufficiently large $k$ and for small enough $\eta_n, \mu_n > 0$ 
\bean
H_{nk}(z)&\leq& \psi(t_{n}, x_{n}+x_{n}z)-\tilde \psi^{\delta}(t'_{n}, y_{n}+y_{n}z)-\psi(t_{n}, x_{n})+\tilde \psi^{\delta}(t'_{n}, y_{n})\\
&&-  nz(x_n - y_n)^2 + \eta_n\\
&\leq& \Delta_n(t_n, x_n + x_nz, y_n + y_nz) - \Delta_n(t_n, x_n, y_n) + \frac12 nz^2 |x_n - y_n|^2 + \eta_n.
\eean
Similarly, 
\bean
G_{nk}(z)&\leq& \psi(t_{nk}, x_{nk}+z)-\tilde \psi^{\delta}(t'_{nk}, y_{nk}+z)-\psi(t_{nk}, x_{nk})+\tilde \psi^{\delta}(t'_{nk}, y_{nk})  + \mu_n \\
&\leq& \Delta_n(t_n, x_n + x_nz, y_n + y_nz) - \Delta_n(t_n, x_n, y_n)  + \mu_n.\\
\eean 
The first terms in both inequalities are dominated by the second terms. The third term in the first inequality is integrable and converges to zero as $n \to \infty$. Hence, 
\bean
&&\limsup_n \frac12 n|x_n-y_n|^2 \int z^2 I_{\{|z | \le 1\}}\Pi(dz) \le 0.
\eean
The last terms $\eta_n$ and $\mu_n$ can be chosen simultaneously with $\delta_n$ to gurantee that for sufficiently large $n$
\bean
&&\mu_n \int I_{\{\delta_n \le |z | \le 1\}}\Pi(dz) \le \frac1n, \quad \eta_n \int I_{\{\delta_n \le |z | \le 1\}}\Pi_{P}(dz) \le \frac1n.
\eean
Hence, applying the dominated convergence theorem,  
\bea
&&\limsup_n\limsup_k \int H_{nk}(z)I_{\{\vert z \vert \le 1\}}\Pi(dz)\le 0,
 \\
&& \limsup_n \limsup_k \int G_{nk}(z)I_{\{\vert z \vert \le 1\}}\Pi_P(dz)\le 0.
\eea 
\fdem 

\medskip
 
As a corollary we get the following uniqueness result. 
 
 \begin{theorem}
\label{unique}  
Let $\Psi$ and $\tilde \Psi$ be two continuous  functions on 
$\bar K:=[0,T]\times \bar \R_+$  which are 
viscosity solutions of (\ref{basiceq}) on $K$  such that  
$$
\Psi(t,0) = \tilde\Psi(t,0), \quad \Psi(t,\infty) =\tilde\Psi(t,\infty), \quad  \Psi(T,u) = \tilde\Psi(T,u).
$$ 
Then $\Psi=\tilde \Psi$.  
\end{theorem}
                               
\acknowledgement  
{The research is funded by the grant of RSF $n^\circ$ 20-68-47030 ``Econometric and probabilistic methods for the analysis of financial markets with complex structure''. }

 \end{document}